\documentclass[12pt]{article}
\usepackage{amsmath, amsfonts, amssymb, amsthm} \parskip1pt
\newtheorem{theorem}{Theorem}[section]
\newtheorem{proposition}[theorem]{Proposition}

\newtheorem{remark}[theorem]{Remark}

\newtheorem{corollary}[theorem]{Corollary}
\newtheorem{lemma}[theorem]{Lemma}

\def\R{I\!\!R} \def\Z{Z\!\!\!Z} \def\T{T\!\!\!T}

\def\tr{{\sf t}}

\begin{document}  
\title{Some two-step and three-step nilpotent Lie groups with small automorphism groups}
\author{S.G. Dani}
\maketitle

Let $G$  be a connected Lie group and  Aut$\,(G)$ be the group of all 
(continuous) automorphisms of $G$. When does 
the action of Aut$\,(G)$ on $G$ have a dense orbit?  It turns out that
if this holds then  $G$ is a nilpotent Lie group; (see~\cite{D2}). However
it does not hold for all nilpotent Lie groups, and a characterization of 
the class of groups for which it holds seems to be a remote possibility.
When $G$ is a vector space then Aut$\,(G)$ is its general linear group, 
and the action has an open dense orbit, namely the complement of 
the zero.  On the other hand there exist
3-step simply connected nilpotent Lie groups such that every orbit of 
the action is a proper closed subset; this holds for the simply connected 
Lie group corresponding to the Lie algebra described in~\cite{DiL}; see 
also $\S\,$\ref{nilauto} below. Among the simply connected groups this 
brings us to considering the question for 
2-step nilpotent Lie groups.  In this note we construct an example of
a 2-step simply connected nilpotent Lie group $G$ for which the action 
of the automorphism group has no dense orbits; see Theorem~\ref{thm}.

We discuss the question also for Lie groups which are not simply
connected, and show that while for all connected abelian groups other than
the circle the automorphism group action has dense orbits, among the 
quotients of the group $G$ as above by discrete central subgroups there are  
nilpotent Lie groups $G'$ such that every orbit of Aut$\,(G')$ on $G'$ 
consists of either one or two cosets of the commutator subgroup $[G',G']$, the 
latter being a closed subgroup;  see Corollary~\ref{corq}.

The method involved also enables us to give an example of a 3-step 
simply connected
Lie group whose automorphism group is nilpotent (its action on the Lie algebra
is by unipotent transformations); see
$\S\,$\ref{nilauto}; the condition implies in particular that all orbits of its
action on $G$ are closed. 
An example of a 6-step simply connected nilpotent Lie group with this property
was given earlier in \cite{Dy}.

Another motivation for studying the automorphism groups of nilpotent
Lie groups comes from the question of understanding which compact
nilmanifolds support Anosov diffeomorphisms and which do not. It is known 
that if $G$ is a  free $k$-step
nilpotent Lie group over a $n$-dimensional vector space $V$ with  $k<n$, 
then for any lattice $\Gamma$ in $G$ the nilmanifold $G/\Gamma$ admits 
Anosov automorphisms; see
\cite{D1}, pp. 558; see \cite{De} for another approach to
the question. There are also other classes of compact nilmanifolds for which 
this holds; see~\cite{AuS}, \cite{D1}, \cite{De}, \cite{DeM}. However 
in general compact nilmanifolds may not admit Anosov
automorphisms. Clearly, if $G$ is a 2-step simply connected nilpotent
Lie group such that $[G,G]$ is one-dimensional then  $G/\Gamma$ can not
admit an Anosov automorphism, for any lattice $\Gamma$. Only a few
other examples of 2-step simply connected Lie groups are known with
this property; see \cite{AuS}, \cite{De}, \cite{M}. We
shall show that for the 2-step nilpotent Lie group $G$ in our example, if
$\Gamma$ is a lattice in $G$ then the nilmanifold $G/\Gamma$ has no
Anosov automorphism, and also no ergodic automorphism; see
Corollary~\ref{coran}.

\section{Lie algebras and Lie groups}\label{Lie}

We recall that a Lie algebra $\cal G$ is said to be nilpotent if  the central 
series $\{ {\cal G}_i\} $ defined by
${\cal G}_0=\cal G$ and ${\cal G}_{i+1} =[{\cal G}, {\cal G}_i]$ for
all $i\geq 0$, terminates, namely there exists $k\geq 1$ such that
${\cal G}_k=0$; if $k$ is the smallest integer for which this holds
then $\cal G$ is said to be $k$-step nilpotent. We say that a
connected Lie group $G$ is a $k$-step nilpotent Lie group if the Lie
algebra of $G$ is a $k$-step nilpotent Lie algebra.

For a Lie group $G$ we shall denote by Aut$\,(G)$ the group of all
continuous automorphisms of $G$, and similarly for a Lie algebra
$\cal G$ we denote by Aut$\,({\cal G})$ the group of all Lie
automorphisms of $\cal G$. 

Let $G$ be a connected nilpotent Lie group. Let $G^{(1)}= \overline{[ G,  G]}$ 
and $Z(G)$ denote the center of $G$. Suppose that $Z(G)$ is contained in 
$G^{(1)}$;  this 
condition holds for the nilpotent groups discussed in the following sections. 
Then for any (continuous) homomorphism $\psi : G/G^{(1)}
\to Z(G) $ the map $\tau:G\to G$ defined by $\tau (g)=g\psi (gG^{(1)})$ 
for all $g\in G$ is an automorphism of $G$; we shall call an automorphism 
arising in this way a {\it shear automorphism}. The class of shear 
automorphisms forms a normal subgroup of Aut$\,(G)$. Every automorphism 
$\tau$ of $ G$ factors to an automorphism of $G/G^{(1)}$; we 
denote the factor of $\tau$ on $G/G^{(1)}$ by $\overline \tau$. Clearly 
$\tau \mapsto \overline \tau$ is a homomorphism of Aut$\,(G)$ into 
Aut$\,(G/G^{(1)})$. We denote by ${\cal A}(G)$ the image of the homomorphism, 
namely, 
$${\cal A}(G)=\{\overline \tau \in \hbox{Aut}\,(G/G^{(1)}) \mid \tau \in 
\hbox{Aut}\,( G)\}.$$ We note that for all shear automorphisms 
as well as all inner automorphisms of $G$ the factor on $G/G^{(1)}$ is 
trivial. 

Suppose now that  $G$ is a simply connected nilpotent Lie group and let 
$\cal G$ be the Lie algebra of $G$. Then $G^{(1)}=[G,G]$, namely $[G,G]$ is 
closed, and $G/[G,G]$ and $[G,G]$ are (topologically isomorphic to) vector 
spaces; they may be identified canonically 
with ${\cal G}/[{\cal G}, {\cal G}] $ and $[{\cal G}, {\cal G}] $ 
respectively. 
Also in this case Aut$\,(G)$ can be realised as Aut$\,({\cal G})$ identifying  
each automorphism with its derivative on $\cal G$. We note also the following:

\begin{lemma}\label{lem:ss}
Let $G$ be a simply connected nilpotent Lie group and $\cal G$ be the Lie
algebra of $G$. Let $\tau \in \hbox{Aut}\,({\cal G})$ and suppose that 
$\overline \tau$ is unipotent (as an element of $GL({\cal G}/[{\cal G}, 
{\cal G}])$). Then $\tau$ is unipotent (as an element of $GL({\cal G})$).
\end{lemma}

\noindent{\it Proof}: It can be seen that the largest $\tau$-invariant 
subspace, say ${\cal G}_0$, of $\cal G$ on which the restriction of $\tau$ 
is unipotent (namely the 
generalized eigenspace for~1 as the eigenvalue) is a Lie subalgebra of $\cal G$. 
Since $\overline \tau$ 
is unipotent ${\cal G}_0+ [{\cal G}, {\cal G}] = \cal G$. Substituting from 
the equation for $\cal G$
on the left hand side successively and using that  $\cal G$ is nilpotent
we deduce that ${\cal G}_0=\cal G$. Therefore $\tau$ is unipotent.   

\medskip
Now let $V$ be a (finite dimensional) vector space and let $W$ be a 
subspace of $\wedge^2\,V$, the second exterior power of $V$. 
We can associate to these canonically a 2-step nilpotent Lie 
group as follows. Let $V'=(\wedge^2\,V)/W$ and ${\cal G}=V\oplus V'$.
We set $[v_1,v_2]=v_1\wedge v_2 \,\hbox{mod}\, W$ for all 
$v_1,v_2 \in V$, and $[x,y]=0$ for all $x\in {\cal G}$ and $y\in V'$.
These relations extend uniquely to a Lie bracket operation on $\cal G$.
Furthermore $\cal G$ is a 2-step nilpotent Lie algebra, with 
${\cal G}/[{\cal G}, {\cal G}] =V$. Conversely every 2-step nilpotent Lie 
algebra $\cal G$ can be realised as a Lie algebra associated to a datum
as above, with $V={\cal G}/[{\cal G}, {\cal G}] $, and $W$ a suitable 
subspace of $\wedge^2\,V$. 

Let $G$ be the simply connected Lie group corresponding to the Lie algebra
$\cal G$ associated to a pair $V,W$ as above. Then $G/[G,G]$ may be
realised canonically as $V$; this identifies the subgroup ${\cal A}(G)$ with 
a subgroup of $GL(V)$. It is easy to see that $x\in GL(V)$ is of the form 
$\overline \tau$ for some $\tau \in \hbox{Aut}\,(G)$ if and only if the 
subspace $W$ of $\wedge^2\,V$ is invariant under the action induced by 
$x$ on $\wedge^2\,V$. It follows in particular that ${\cal A}(G)$ is an 
algebraic subgroup of $GL(V)$. 

\medskip
We note also the following:

\begin{proposition}\label{prop}
Let $G$ be a simply connected nilpotent Lie group. Then every orbit of 
Aut$\,(G)$ on $G$  is open in its closure. In particular every 
dense orbit is open. Consequently, if the ${\cal A}(G)$-action on $V$ has 
no open orbit then the Aut$\,(G)$-action on $G$ has no dense orbit. 
\end{proposition}

\noindent{\it Proof}: Let $\cal G$ be the Lie algebra of $G$. Then
Aut$\,({\cal G})$ is an algebraic subgroup of
$GL(\cal G)$ and this implies that the orbits of its action on $\cal G$
are open in their closures (see \cite{BR}). Since $G$ is simply connected
and nilpotent the exponential map is a diffeomorphism of $\cal G$ to $G$
and it is equivariant under the actions of the respective automorphism groups. 
Therefore for the action of Aut$\,(G)$ on $G$ also every orbit is open 
in its closure. The second assertion is immediate from the first. The last 
assertion follows from the fact that the ${\cal A}(G)$-action on $V$ 
is a factor of the Aut$\,(G)$-action on $G$. 

\section{A representation}\label{rep}

In this section we describe a construction of a representation and 
prove some properties. The results will be used in the later sections
to give examples of nilpotent Lie algebras. 

Let $\cal S$ be the vector space of $3\times 3 $ symmetric matrices
with real entries. For each $k,l=1,2,3$ we denote by $E_{kl}$  the
$3\times 3$ matrix in which the $(k,l)$-entry (in the $k$th row and
$l$th column) is $1$ and all other entries are $0$.  Let $V$ be the
subspace of $\cal S$ defined by $$V= \{\Sigma\, \sigma_{kl}E_{kl}
\in {\cal S}\mid
\sigma_{22}=\sigma_{13}+\sigma_{31}=2\sigma_{13}\}.$$ Let $\sigma_1,
\dots , \sigma_5 $ be the elements defined as follows, forming a
basis of $V$: $$\sigma_1 =2E_{11}, \ \sigma_2=E_{12}+E_{21}, \
\sigma_3=E_{13}+2E_{22}+E_{31}, \ \sigma_4 =E_{23}+E_{32}, \ \hbox{
and } \sigma_5= 2E_{33}.$$
Let $W$ be  the subspace of $\wedge^2\,V$ spanned by the  three elements
$\sigma_1\wedge \sigma_4 -\sigma_2\wedge\sigma_3$, 
$ \sigma_1\wedge\sigma_5 - \sigma_2\wedge
\sigma_4$ and $\sigma_2\wedge\sigma_5 - \sigma_3\wedge \sigma_4$. 

Let $\delta$, $\nu^+$ and $\nu^-$ be the matrices defined by $$\delta
=\begin{pmatrix} 2 & 0 & 0 \\ 0 & 0 & 0 \\ 0 & 0 & -2\\
\end{pmatrix}, \ \nu^+ =\begin{pmatrix} 0 & 1 & 0 \\ 0 & 0 & 1 \\ 0 &
0 & 0 \\\end{pmatrix},\hbox{ and } \nu^- =\begin{pmatrix}0 & 0 & 0 \\
1 & 0 & 0 \\ 0 & 1 & 0 \\\end{pmatrix}.$$ \\ 
The space spanned by the three matrices is a Lie subalgebra of the
Lie algebra of $3\times 3$ matrices of trace $0$, namely the Lie
algebra of $SL(3,\R)$. We denote the Lie subalgebra by $\cal H$. We
note that it is isomorphic to the Lie algebra of $SL(2,\R)$, with 
$\delta$, $\nu^+$ and $\nu^-$ corresponding to the standard basis of 
the latter; in particular
the subspace spanned by $\delta$ is a Cartan subalgebra. Let $H$ be
the connected Lie subgroup of $SL(3,\R)$ corresponding to $\cal H$.
Consider the action of $H$ on $\R^3$ given  by restriction of
the natural action of $SL(3,\R)$ on $\R^3$. We see that there is no
proper nonzero subspace on $\R^3$ invariant under both $\exp \nu^+$ and
$\exp \nu^-$. It follows therefore that the action of $H$ on $\R^3$ is 
irreducible.  Since
$H$ contains an element, viz $\exp\,\delta$, with distinct eigenvalues
this implies that every element of the center of $H$ acts by scalar
multiplication by real numbers, and considering the determinant we
see that it must be trivial. Thus $H$ has trivial center. This shows
that $H$ is Lie isomorphic to $PSL(2,\R)$, the adjoint group of
$SL(2,\R)$.

Consider the action of $SL(3,\R)$ on $\cal S$ given by $(x,\sigma)
\mapsto x \sigma x^\tr$ for all $x\in SL(3,\R) $ and $\sigma \in \cal
S$, where $x^\tr$ denotes the transpose of $x$. It is straightforward
to verify that the subspace $V$ is invariant under the action  of
$H$ (obtained by restriction); it suffices to verify that $V$ is invariant 
under the corresponding action of the Lie algebra, given by $(\xi, \sigma) 
\mapsto \xi \sigma + \sigma \xi^\tr$, for all $\xi \in \cal H$ and 
$\sigma \in \Sigma$, and furthermore it is enough to consider $\xi 
= \delta, \nu^+$ and $\nu^-$. The elements $\sigma_1, \dots
,\sigma_5$ are weight vectors with respect to $\delta$ with weights
$4, 2, 0, -2$ and $-4$ respectively, and the latter being distinct
shows that the representation of $H$ over $V$ is irreducible.

We shall realize $H$ as a subgroup of $GL(V)$ by identifying each $h$
in $H$ with its action on $V$ (this is indeed an injective
correspondence). We denote by $D$ the one-parameter subgroup of
$GL(V)$ corresponding to the one-parameter subgroup $\{\exp t\delta
\mid t\in \R\}$.  We note that for any nontrivial element of $D$,
$\sigma_1, \dots, \sigma_5$ are eigenvectors with distinct eigenvalues,
and hence the centralise of $D$ in $GL(V)$ consists of diagonal 
matrices with respect to the basis $\sigma_1, \dots, \sigma_5$. 

Let ${\cal A}$ be the subgroup of consisting of all $x$ in $GL(V)$ such 
that $W$ is invariant under the action induced by $x$ on $\wedge^2\,V$. 
A straightforward  computation shows that   $H$ (viewed as a subgroup 
of $GL(V)$) is contained in ${\cal A}$. We note that $\cal A$ is
an algebraic subgroup of $GL(V)$.  Let ${\cal A}_u$ be the 
unipotent radical of $\cal A$, namely the largest normal subgroup consisting 
of unipotent elements. Then the  set
of common fixed points of the action of ${\cal A}_u$ on $V$ is a
nonzero subspace invariant under $H$, and since  the $H$-action on
$V$ is irreducible it follows that ${\cal A}_u$ fixes all points of
$V$, which means that ${\cal A}_u$ is trivial. Therefore $\cal A$ is 
reductive.

\begin{theorem}\label{thm}
${\cal A}=ZH$, where $Z$ is the subgroup of $GL(V)$ consisting  
of scalar multiplications by real numbers. The ${\cal
A}$-action on $V$ has no open orbit.  
\end{theorem}

\noindent{\it Proof}: Let ${\cal A}^0$ be the connected component of 
the identity in $\cal A$ and
let $S=[{\cal A}^0,{\cal A}^0]$, the commutator subgroup. As ${\cal
A}^0$ is reductive it follows that  $S$ is a semisimple subgroup of
$GL(V)$; (see \cite{H}, \cite{V}). Since $H$ is a simple Lie subgroup of ${\cal
A}^0$ it follows that it is contained in $S$. We shall show that
$S=H$. Firstly let $C$ be the maximal compact connected normal
subgroup of $S$. Since $S$ is semisimple, there exists a closed
connected normal subgroup $S'$ such that $S=S'C$ and $S'\cap C$ is a
finite subgroup contained in the center of $S$; both $S'$ and $C$
being connected subgroups this implies in particular that every
element of $S'$ commutes with every element of $C$.  Since $H$ is a
noncompact connected simple Lie group it has no nontrivial
homomorphism into a compact Lie group. Applying this to the quotient
homomorphism of $S$ onto $S/S'=C/(C\cap S')$ we conclude that $H$ is contained 
in $S'$. Therefore $C$ is contained in the centralizer of $D$. Since the
centralizer  is diagonalizable and $C$ is a compact connected subgroup it 
follows that $C$ is trivial. Thus $S$ has no compact normal subgroups 
of positive dimension.

Now recall that the subspace $W$ of $\wedge^2\,V$ is invariant under
the action of $\cal A$ and in particular that of $S$, and let
$\varphi:S\to GL(W)$ be the representation of $S$ induced by the
action. Let $K$ be the connected component of the identity in the
kernel of $\varphi$.  Then $K$ is a connected normal subgroup of $S$.
Let $A$ be a maximal subgroup of $S$ such that $D$ is contained in it
and its adjoint action  on the Lie algebra of $S$ is diagonalisable
over $\R$.  Then $A$ intersects (by Aswan decomposition) any
noncompact connected normal subgroup nontrivially, and since $S$ has
no compact normal subgroups it follows that $A\cap K$ is nontrivial,
unless $K$ is trivial. Let $\alpha\in A\cap K$. Then  $\alpha$ commutes 
with all elements of $D$,
and hence $\sigma_1, \dots ,\sigma_5$ are eigenvectors of
$\alpha$; let $\lambda_1, \dots , \lambda_5$ be the corresponding
eigenvalues. Since $\alpha \in K$ the action induced by $\alpha$ on
$W$ is trivial. This means that the vectors $\sigma_1\wedge \sigma_4
-\sigma_2\wedge\sigma_3$, $\sigma_1\wedge\sigma_5 - \sigma_2\wedge \sigma_4$ 
and $\sigma_2\wedge\sigma_5 - \sigma_3\wedge \sigma_4$ are all 
fixed under the action of $\alpha$.
Each $\sigma_i\wedge \sigma_j$ is an eigenvector of the action of
$\alpha$ on $\wedge^2\,V$ and hence the three vectors as above being
fixed  implies that the vectors $\sigma_1\wedge \sigma_4$, 
$\sigma_1\wedge\sigma_5$, $\sigma_2\wedge\sigma_3$, $\sigma_2\wedge \sigma_4$,
$\sigma_2\wedge\sigma_5$ and $\sigma_3\wedge \sigma_4$ 
 are fixed individually.  Therefore for the
eigenvalues we get the relations $\lambda_1\lambda_4=\lambda_1\lambda_5=
\lambda_2\lambda_3=\lambda_2\lambda_4=\lambda_2\lambda_5=\lambda_3\lambda_4
=1$, which is possible only if $\lambda_i=1$ for
all $i=1,\dots ,5$. Since $\alpha$ is diagonalizable this shows that
$\alpha$ is the identity. Thus $A\cap K$ is trivial and
hence, as noted above, $K$ is trivial. Therefore the kernel of
$\varphi$ is discrete.

Clearly  $\varphi (S)$ is a connected semisimple Lie subgroup of
$GL(W)$ containing $\varphi (H)$. Since $W$ is 3-dimensional this
implies that $\varphi (S)$ is either $\varphi (H)$ or $SL(W)$;
(inspection of the weights for the action of $\varphi(D)$ on the Lie
algebra of $SL(W)$, via the adjoint representation of $SL(W)$, shows that the
$\varphi (H)$-action has only two irreducible components, one of them
being the Lie subalgebra of $\varphi (H)$ itself;  this means that
there is not even a {\it subspace} invariant under the $\varphi
(H)$-action lying strictly between the Lie algebras of $\varphi (H)$
and $SL(W)$, which in particular implies the assertion here; the
assertion can also be proved in other ways). Since
$\varphi$ has discrete kernel this shows that the Lie algebra of $S$ is 
isomorphic to that of either $H$ or $SL(3,\R)$. We note that the
Lie algebra of $SL(3,\R)$ has no irreducible 5-dimensional
representation; this can be deduced from the corresponding statement
over the field of complex numbers, the latter being easy to see from
the classification theory of representations of semisimple Lie
algebras; see \cite{S}. Recall that the $H$-action on $V$ is
irreducible. Since $S$ contains $H$, its action on $V$ is
irreducible, and since  $V$ is 5-dimensional the preceding
observation implies that the Lie algebra of $S$ can not be isomorphic to
that of $SL(3,\R)$. Therefore from the alternatives above we now get that
$S$ is locally isomorphic to $H$. Since $S$ is connected and $H$ is a
subgroup of $S$ this implies that $S=H$.

Since $H=S=[{\cal A}^0,{\cal A}^0]$, in particular $H$ is normal in 
${\cal A}$. Now let $\alpha \in \cal A$ be arbitrary and consider the 
automorphism of $H$ induced by the conjugation action of $\alpha$. 
Since $H$ is Lie isomorphic to $PSL(2,\R)$ every automorphism of 
$H$ is inner. Hence there exists $h_0 \in H$ such that 
$\alpha h \alpha^{-1}=h_0hh_0^{-1}$ for all $h\in H$. Then
$h_0^{-1}\alpha$ commutes with every element of $H$ and
particular with all elements
of $D$. Therefore $h_0^{-1}\alpha$ is diagonalisable and in particular
all its  eigenvalues  are real.
Since the action of $H$ on $V$ is irreducible it now follows that
$h_0^{-1}\alpha$ acts by scalar multiplication by a (nonzero) real 
number. Thus $h_0^{-1}\alpha\in  Z$ in the notation as in the
hypothesis.  Then $\alpha \in HZ=ZH$, and since $\alpha \in \cal A$
was arbitrary we get that ${\cal A}$ is contained in $ZH$. On the other
hand $H$ and $Z$ are contained in $\cal A$ and therefore we have 
${\cal A}=ZH$. This proves the first assertion in the theorem. The 
second assertion follows from the fact that $ZH$ is 4-dimensional 
while $V$ is 5-dimensional. 

\medskip
We note also the following simple fact about representations of $PSL(2,\R)$:

\begin{lemma}\label{lem}
Let $\rho :PSL(2,\R) \to GL(Q)$ be a representation 
of $PSL(2,\R)$ over a $\R$-vector space $Q$ of dimension $n \geq 5$,
such  that no nonzero point of $Q$ is  fixed by the image of $\rho$. Let $E$ 
be the set of points $p$ in $Q$  such that $p$ is an eigenvector of $\rho (g)$ 
for some nontrivial  $g \in PSL(2,\R)$. Then $E$ is a union of countably many 
smooth submanifolds of dimension at most $(n+4)/2$. In particular if $L$ is 
a subspace of $Q$ of dimension $m> (n+4)/2$ and $\lambda$ is the Lebesgue 
measure on $L$ then $\lambda (E\cap L) =0$. 
\end{lemma}

\noindent{\it Proof}: We realize $PSL(2,\R)$ as the group $H$ as above
(this is just for notational convenience). 
For each $t\in \R$ let $d_t=\exp t\delta$ and $u_t=\exp t\nu^+$, where 
$\delta $ and $\nu^+$ are the matrices as defined earlier. Let $\{k_t\}$ 
be a periodic one-parameter subgroup of $H$, say with period $2\pi$. Then 
every element of $H$ is conjugate to one of $d_t$, $u_t$ or $k_t$ for some 
$t\in \R$. 
Let $R$ be the countable subset of $H$ consisting of $d_1, u_1, k_1$ and $k_r$ 
for all $r$ of the form $2\pi /m$ with $m$ a positive integer. 

For any $h \in H$ let $E(h)$ denote 
the set of points in $L$ which are eigenvectors of $\rho (h)$. Then $E(h)$ is a 
union of finitely many vector subspaces of $L$; we shall denote by $d(h)$
the maximum of the dimension of the subspaces contained in $E(h)$. 
Now consider any nontrivial element $h$ in $H$. We can find a 
conjugate $h'$ of $h$ such that the subgroup generated by $h$ and $h'$ is
dense in $H=PSL(2,\R)$. Then every point of $E(h)\cap E(h')$ is an eigenvector 
of all elements of $\rho (H)$, and since the latter is a simple Lie group 
the points must be fixed by $\rho (H)$. The condition in the hypothesis 
therefore implies that $E(h)\cap E(h')=0$. Since $h$ and $h'$ are conjugates
$d(h)=d(h')$ and hence the preceding conclusion implies that 
$d(h) \leq n/2$. 

Now let $p\in E$ and $h\in H$ be such that $p$ is an eigenvector of $\rho (h)$. 
Suppose first that $h$ is of infinite order. There exists $t\in \R$ and $g\in H$  
such that $ghg^{-1}$ is one of $d_t$, $u_t$ or $k_t$.  Then $\rho (g)(p)$ is an
eigenvector of $d_t$, $u_t$ or $k_t$. This implies that it is an 
eigenvector of one of the elements $d_1$, $u_1$ or $k_1$ respectively; in 
respect of the last alternative note that since $k_t$ is of infinite order the 
closure of the cyclic subgroup generated by it contains $k_1$.
Thus there exists $x\in R$ such that $\rho (g)(p) \in E(x)$, and so 
$p \in \rho (H)E(x)$. 
Now suppose that $h$ as above is of finite order. Then $h$ is conjugate to $k_r$
for some $r$ of the form $2\pi/m$ for a positive integer $m$. Therefore in 
this case also we get that  $\rho (g)(p) \in E(x)$ for some $x\in R$. Thus 
$E=\cup_{x\in R}\,\rho (H)(E(x))$. Each $E(x)$ is a finite union of vector 
subspaces which are  invariant under 
a one-parameter subgroup of $H$ and have dimension at most $n/2$. It follows 
that each $\rho (H)(E(x))$ is a finite union smooth 
submanifolds of dimension at most $(n/2)+2=(n+4)/2$.  This proves the Lemma. 

\medskip
In respect of the conclusion as in the lemma the situation in  low dimensions
is as follows. For dimensions 2 and 4 there are no representations of $PSL(2,\R)$ 
with no nonzero common fixed points. In dimension~3 the realization of $H$ (as above)
as a subgroup of $SL(3,\R)$ is (up to equivalence) the only representation with no 
common fixed points. For this representation the set of points which are 
eigenvectors of nontrivial elements form a cone in $\R^3$ with nonempty interior. 

\section{A 2-step nilpotent Lie group}\label{2step}

We now use Theorem~\ref{thm} to give an example of a simply connected 
2-step nilpotent Lie group $G$ such that the action of Aut$\,(G)$ on $G$ has 
no dense orbit.

Let $V$ and $W \subset \wedge^2\,V$ be as in the preceding section and 
let $\cal G$  be the 2-step nilpotent Lie algebra associated with the pair $V$, $W$ as 
in $\S\,$\ref{Lie};  namely ${\cal G}=V\oplus V'$, where  
$V'=(\wedge^2\,V)/W$ and the Lie bracket operation is determined by the
conditions $[v_1,v_2] =(v_1\wedge v_2)$ (mod $W$) for all $v_1,v_2 \in V$,
and  $[x,y]=0$ for $x\in \cal G$ and $y\in V'$.
Let $G$ be the simply connected nilpotent Lie group corresponding to  $\cal G$. 
Then we have

\begin{corollary}\label{cormain}
The action of  Aut$\,(G)$ on $G$ has no dense orbit. 
\end{corollary}

\noindent{\it Proof}: For the Lie algebra $\cal G$ as above, corresponding to $G$, 
the subgroup ${\cal A}({\cal G})$ as defined in $\S\,$\ref{Lie} is the same as 
the subgroup $\cal A$ as in Theorem~\ref{thm}.  Since by Theorem~\ref{thm}
the action of $\cal A$ on 
$V$ has no open orbit, by Proposition~\ref{prop} the Aut$\,(G)$-action on $G$ 
has no dense orbit.

\medskip
Using Theorem~\ref{thm} one can also describe Aut$\,(G)$ as follows. 
To each $h \in H$ corresponds a Lie automorphism of $\cal G$ such that 
the restrictions to $V$ and $V'=(\wedge^2\,V)/W$ are respectively given 
by the $h$-action on $V$ and the factor of the $h$-action on $\wedge^2\,V$.
We denote by $\tilde H$ the group of automorphisms of $G$ 
for which the associated Lie automorphism corresponds to some $h\in H$ 
as above. Then $\tilde H$ is a Lie subgroup of Aut$\,(G)$, canonically Lie
isomorphic $H$. Similarly each $z\in Z$ corresponds canonically to an 
automorphism of $G$, and we shall denote by $\tilde Z$ the subgroup of 
Aut$\,(G)$ corresponding to $Z$. Also let $\Psi$ be the group of shear 
automorphisms of $G$ (see $\S\,$\ref{Lie}). 

\begin{corollary}\label{coraut}
Aut$\,(G)=(\tilde Z \tilde H)\Psi$, semidirect product 
(with $\Psi$ as the normal subgroup).
\end{corollary}

\smallskip
\noindent{\it Proof}: We note that  $(\tilde Z \tilde H)\Psi$ is a subgroup of Aut$\,(G)$.
Now consider the homomorphism $\tau \mapsto \overline \tau$ of Aut$\,(G)$ into 
$GL(V)$. It is easy to see that its kernel is $\Psi$, and its image is $\cal A$ which
by  Theorem~\ref{thm} is $ZH$.  As the subgroup $(\tilde Z \tilde H)\Psi$  of Aut$\,(G)$ 
contains $\Psi$, and its image under the homomorphism as above is $ZH$, this shows 
that it must be the whole.  Since the map $\tau \mapsto \overline \tau$
is injective on $\tilde Z \tilde H$ it follows that $(\tilde Z \tilde H)\cap \Psi$ is trivial. 
This proves the corollary. 

\begin{remark} {\rm The action of $H$ on $\wedge^2\,V$
(notation as before) decomposes into  
two irreducible components, one being on the subspace $W$ as above
and another on the subspace, say $W'$, spanned by the $H$-orbit of
$\sigma_1\wedge \sigma_2$; this can be deduced by inspection of the
set of weights. Let ${\cal G}'$ be the 2-step nilpotent Lie algebra 
associated to the pair $V$ and $W'$, as described in $\S\,$\ref{Lie}. Then 
with some modifications in the above argument it can be
shown that Aut$\,( G')$ has no open orbit on $G'$; in
this respect we note mainly that  $W$ is also invariant, 
along with $W'$,  under the action of 
${\cal A}(G')$ on $\wedge^2\,V$. This yields an
$8$-dimensional example (as against $12$ of $\cal G$) for which the
automorphism group action has no dense orbit. The example as in 
Theorem~\ref{thm} on the other hand admits a simpler presentation.}
\end{remark}

\section{Lie groups with nilpotent automorphism groups}\label{nilauto}

An example of a 6-step simply connected
nilpotent Lie group such that Aut$\,(G)$ is nilpotent was given in \cite{Dy}. 
Earlier in~\cite{DiL} an example of a 3-step nilpotent 
Lie algebra $\cal G$ was constructed for which  all derivations are 
nilpotent; the condition implies that Aut$\,({\cal G})^0$,
the connected component of the identity in  Aut$\,({\cal G})$, consists of
unipotent elements (viewed as elements in $GL({\cal G})$), and hence is a 
nilpotent Lie group; 
Aut$\,({\cal G})$ itself is however not nilpotent for that example, as 
has been remarked in \cite{Dy}.
In this section we shall describe a class of examples of simply connected 
3-step nilpotent Lie 
groups, related  to the results in the preceding sections, for which the all 
automorphisms are unipotent,  and the automorphism groups are nilpotent.
The following implication of this to orbits of the automorphism group 
action may be borne in mind. 

\begin{remark}\label{rem:co}
{\rm Let $\cal G$ be a nilpotent Lie algebra such that  Aut$\,({\cal G})^0$ 
consists of unipotent elements. Then  Aut$\,({\cal G})^0$ is a unipotent
algebraic subgroup of $GL({\cal G})$ and hence all its orbits on $\cal G$
are closed (see~\cite{H}). Since  Aut$\,({\cal G})$ is an algebraic group
it has only finitely many connected components and hence we get furthermore
that all orbits of  Aut$\,({\cal G})$ on $\cal G$ are closed. As the 
Aut$\,(G)$-action on $G$ and the  Aut$\,({\cal G})$-action on $\cal G$ are
topologically equivalent (via the exponential map) it follows that all 
orbits of Aut$\,(G)$ on $G$ are closed. 
Thus in this case we have a situation stronger than the
orbits not being dense, that was obtained in $\S\,$\ref{2step}.} 
\end{remark}

We note also the following, which shows that the 3-step condition in the
above assertions is optimal. 

\begin{remark}\label{rem:2s}
{\rm Let $G$ be a 2-step simply connected nilpotent Lie group. Then for 
the action of Aut$\,(G)$ on $G$ the closure of every orbit contains the
identity element; in particular orbits other than that of the identity
element are not closed; (in 
this case the automorphism group, or even its connected component of 
the identity, does not consist entirely of unipotent elements (viewed
as automorphisms of the Lie algebra). The assertion is immediate from 
the fact that all scalar transformations of $V=G/[G,G]$ belong to 
${\cal A}(G)$; this shows that the closure of every orbit of   
Aut$\,({\cal G})$ on $\cal G$, where $\cal G$ is the Lie algebra of $G$, 
contains the zero element of $\cal G$; this is equivalent to the 
statement as above.} 
\end{remark}

We now proceed to describe the 3-step nilpotent Lie algebras. We shall 
follow the notation as in $\S\,$\ref{rep}. 
For $1\leq i <j\leq 5$ let $p_{ij} =\sigma_i\wedge \sigma_j \,\hbox{mod}\,(W)$.
We note that $p_{14}=p_{23}$, $p_{15}=p_{24}$, and $p_{25}=p_{34}$, and 
that $\{p_{12}, p_{13}, p_{14}, p_{15}, p_{25}, p_{35}, p_{45}\}$ 
is a basis of $V'$. Let $L$ be the subspace of $V'$ spanned 
by $\{p_{13}, p_{14}, p_{15}, p_{25}, p_{35}, p_{45}\}$. Then 
by Lemma~\ref{lem} there exists (a nonzero element) $p\in L$ which is not an 
eigenvector of the action of any nontrivial element of $H$ (under the 
$H$-action on $V'$).  Now we set ${\cal N}=V\oplus V'$
and we define a Lie bracket operation by the relations $[\sigma_i,\sigma_j]
=\sigma_i \wedge \sigma_j$ (mod $W$) for all $i,j=1, \dots ,5$,
$[x, x']=0$ for all $x,x'\in V'$, $[\sigma_1,p_{12}]= p$, and $[\sigma_k, p_{ij}]=0$ if 
either $k\neq 1$ or 
$(i,j)\neq (1,2)$; the product defined on the basis elements extends 
uniquely to a Lie bracket operation on $\cal N$ as above; we note in particular
that by definition $[\sigma_i, [\sigma_j,\sigma_k]]=0$ if $i,j$ and $k$ are distinct
indices, from which one can see that the Jacobi identity holds for the 
product. Furthermore, since $p$ is contained in the subspace $L$ as above it
follows that $p$ is contained in the center of $\cal N$.  Since $p$ spans 
$[\cal N, [\cal N,\cal N]]$ this shows that $\cal N$ is a 
3-step nilpotent Lie algebra. Let $N$ be the simply connected (nilpotent) Lie 
group corresponding to $\cal N$. We prove the following:

\begin{corollary}\label{cor3}
The subgroup ${\cal A}(N)$ consists only of the identity element. In
particular 
Aut$\,( N)$ consists of  unipotent elements, when realised (canonically) as a 
subgroup of $GL({\cal N})$.
\end{corollary}

\noindent{\it Proof}: Let 
$\tau \in \hbox{Aut}\,({\cal N})$. Then $\tau$ factors to a linear 
transformation $x=\overline  \tau \in GL(V)$. Since
$[\sigma_i,\sigma_j]
=\sigma_i \wedge \sigma_j$ (mod $W$) for all $i,j=1, \dots ,5$, it follows that
the action of $x$ on $\wedge^2\,V$ leaves invariant the subspace $W$.
Hence by Theorem~\ref{thm} $x \in ZH$. Let $z\in Z$ and $h\in H$ be such that 
$x=zh$. Since the subspace spanned by $p$ equals $[\cal N, [\cal N,\cal N]]$
it is invariant under all automorphisms of $\cal N$. Thus $p$ is an 
eigenvector of $x$ (for the action on $V'$). Since the elements of $Z$ act as
scalars on $V'$, it follows that $p$ is an eigenvector of $h$. Since by choice
$p$ is not an eigenvector of any nontrivial element of $H$, we get that $h$ is
trivial. Thus  $\overline \tau = x=z$, a scalar transformation, say multiplication 
by $\lambda$.  Then the action of $\tau $ on $V'$ is given by multiplication by 
$\lambda^2$. Since $[\sigma_1,p_{12}]= p$, we have $\lambda^2 p=\tau (p) 
=\tau ([\sigma_1 ,p_{12}])=[\tau (\sigma_1),\tau (p_{12})]=\lambda^3 
[\sigma_1, p_{12}]=\lambda^3p$. Since $p$ is a nonzero element this implies that  
$\lambda=1$. This means that $z$ is trivial, and therefore $\overline \tau$ is 
trivial. 
Thus $\overline \tau$ is trivial for all $\tau \in \hbox{Aut}\,({\cal N})$.  
This shows that ${\cal A}(N)$ is trivial. The second assertion follows from  
this, together with Lemma~\ref{lem:ss}. Thus proves the corollary. 

\medskip
We deduce also the following, showing that under a slight further condition
the automorphism group is `minimum possible'.

\begin{corollary}\label{cor:min}
Let $N$ and $\cal N$ be as above. Suppose the point $p$ in the definition of the
Lie algebra structure on $\cal N$ is not contained in the subspace spanned by   
$\{p_{14}, p_{15}, p_{25}, p_{35}, p_{45}\}$. Let $N^*$ denote subgroup 
of Aut$\,(N)$ consisting
of all inner automorphisms and $\Psi$ be the subgroup consisting of all shear 
automorphisms of $N$. Then Aut$\,(N)=N^*\Psi$.
 \end{corollary}

\noindent{\it Proof}: We shall identify Aut$\,(N)$ with Aut$\,({\cal N})$ as 
 before. As an algebraic subgroup consisting of unipotent elements
Aut$\,({\cal N})$ is 
a connected Lie group. The Lie algebra $\cal D$ consisting of all derivations 
of $\cal N$ is the Lie algebra of Aut$\,({\cal N})$.  The subgroups $N^*$ and 
$\Psi$ are  
normal Lie subgroups of  Aut$\,({\cal N})$. Therefore to prove the corollary 
it suffices 
to show that the Lie subalgebras of $\cal D$ corresponding to $N^*$ and $\Psi$ 
span $\cal D$. The Lie subalgebra corresponding to $\Psi$ consists of all 
derivations $\delta$ such that $\delta ({\cal N})$ is contained in $L$. We shall 
show that this Lie subalgebra and  ad$\,\sigma_1$ and ad$\,\sigma_2$ together
span $\cal D$ as a vector space. This would complete the proof. 

Let $\delta$ be any derivation of  $\cal N$. In view of Corollary~\ref{cor3} its
factor on ${\cal N}/[\cal N,\cal N]$ is trivial, and hence $\delta ({\cal N})$ is
contained in $V'$. Hence there exist $a_1,a_2 \in \R$ such that $\delta (\sigma_i)
\in a_ip_{12} +L$,  $i=1, 2$. Then $\delta'=\delta - a_2 
(\hbox{ad}\,\sigma_1) + a_1(\hbox{ad}\,\sigma_2)$ is a derivation such that 
$\delta' (\sigma_1), \delta'(\sigma_2)\in L$. It suffices to show that $\delta'$
 belongs to 
the Lie subalgebra of $\Psi$, namely that $\delta' ({\cal N})$ is contained in $L$. 
Changing notation we shall assume that $\delta (\sigma_1), \delta (\sigma_2)\in L$
and deduce that $\delta ({\cal N}) \subset L$.

As $\delta ({\cal N}) \subset V'$, for $i,j \geq 2$ we  have $\delta (p_{ij})
=\delta ([\sigma_i,\sigma_j])=[\delta \sigma_i,\sigma_j] +[\sigma_i,\delta \sigma_j]
=0$. Since
$p_{14}=p_{23}$ and $p_{15}=p_{24} $, it follows that $\delta (p_{14})=
\delta (p_{15})=0$. We note that $\delta (p_{12})=\delta ([\sigma_1,\sigma_2])
=[\delta \sigma_1, \sigma_2]+[\sigma_1,\delta \sigma_2] \in [{\cal N},[{\cal N},
{\cal N}]]  \subset L$. Therefore $\delta (p)=\delta ([\sigma_1,p_{12}])=
[\delta \sigma_1, p_{12}]+[\sigma_1,\delta p_{12}]=0$.
Let $L'$  be the subspace spanned by 
$\{p_{14}, p_{15}, p_{25}, p_{35}, p_{45}\}$. Then from what we have 
seen $\delta (x)=0$ for all $x\in L'$. By 
the assumption on $p$ there exists $\lambda \neq 0$ such that $p=\lambda p_{13}+x$
for some $x\in L'$. As $\delta (p)=0$ and $\delta (x)=0$, this implies that 
$\delta (p_{13})=0$. Thus we have $[\sigma_1,
\delta (\sigma_j)]=\delta ([\sigma_1,\sigma_j])=\delta (p_{1j})=0$ for $j= 3,4,5$.  
This shows that $\delta (\sigma_j) \in L$ for 
$j=3,4,5$.  Since by assumption $\delta (\sigma_1), 
\delta (\sigma_2)\in L$, we now have $\delta (\sigma_i) \in L$ for all 
$i=1,\dots, 5$. Therefore $\delta ({\cal N})$ is contained in $L$. This 
completes the proof.

\section{Non-simply connected Lie groups}\label{nsc}

In this section we shall discuss the analogous questions for Lie
groups which are not necessarily simply connected. We begin with 
the following observation. 

\begin{proposition}\label{prop:abel}
Let $G$ be a connected abelian Lie group. Then the action of 
Aut$\,(G)$ on $G$ has a dense orbit if and only if $G$ is not 
(topologically isomorphic to) the circle group. 
\end{proposition}

\noindent{\it Proof}: It is well known that for $m\geq 2$ the
$m$-dimensional torus admits automorphisms with dense orbits (see~\cite{W}
for instance). Now let $G=V\times C$, where $V$ is a vector space
of dimension $n\geq 1$, $C$ is the torus of dimension $m\geq 0$. Let $\cal H$
be the set of all continuous homomorphisms of $V$ into $C$. For 
each $\tau \in GL(V)$ and $\psi \in {\cal H}$ we get a
continuous  automorphism
of $G$ defined by $(v,c)\mapsto (\tau (v),c\psi (v))$ for all $v\in V$
and $c\in C$. The automorphisms arising in this way form a subgroup 
of Aut$\,(G)$ (it is in fact the identity component of the latter). 
We recall that there exist homomorphisms $\psi:V\to C$ such that $\psi (V)$
is dense in $C$. Using this 
it is straightforward to verify that under the action of the group of
automorphisms 
as above the orbit of any element of the form $(v,c)$ with $v\neq 0$ is 
dense in $G$. In particular the Aut$\,(G)$-action on $G$ has dense orbits. 
Finally, the circle group has only two
automorphisms and hence in this case the Aut$\,(G)$-action has no dense orbit. 
This proves the proposition.  

For Lie groups $G$ of the form $\R^n\times \T^m$ with $m\neq 1$ there 
exist abelian subgroups of Aut$\,(G)$, with finite rank, whose
action on $G$ has dense orbits; see~\cite{D3} for more precise results
in this respect. 

\medskip
In the context of the  results in the last section and 
Proposition~\ref{prop:abel} one may ask whether there exist  
connected 2-step nilpotent Lie groups $G$ such that all orbits of 
the Aut$\,(G)$-action on $G$ are closed. By Remark~\ref{rem:2s} such a 
group is necessarily non-simply connected. The following corollary shows
that there are groups with this property among quotients of the group $G$ 
as in $\S\,$\ref{2step}, by discrete central subgroups. 

\begin{corollary}\label{corq}
Let $G$ be the simply connected 2-step nilpotent Lie group 
as in Corollary~\ref{cormain}. Then there exists a subset of $E$ of Lebesgue 
measure $0$ in $[G,G]$ such that if $\theta \in [G,G]$, $\theta \notin E$ and 
$\Theta$ is the cyclic subgroup of $G$ generated by $\theta$, then for the 
Lie group
$G'=G/\Theta$ the following holds: for any automorphism $\tau$ of $G'$ the 
factor of $\tau$ on $G'/[G',G']$ is $\pm I$, where $I$ is the identity 
transformation; 
furthermore, orbits of the action of Aut$\,(G')$ on $G'$ consist of either 
one or two cosets of $[G',G']$ in $G'$.
\end{corollary}

\noindent{\it Proof}: Let  $\tilde H$ be the group of automorphisms as 
in $\S\,$\ref{2step} and 
consider the $\tilde H$-action on $[G,G]$. Recall that the latter may be viewed as
a vector space and that the
action is irreducible; in particular there is no nonzero point fixed by the whole 
of $\tilde H$. Since $\tilde H$ is Lie isomorphic
to $PSL(2,\R)$ and $[G,G]$ is 7-dimensional, by Lemma~\ref{lem} there exists 
a subset $E$ 
of Lebesgue measure 0 such that points outside $E$ are not eigenvectors of 
any nontrivial 
element of $\tilde H$. Now let $\theta \in [G,G]$, $\theta \notin E$ and 
$\Theta$ be the
cyclic subgroup generated by $\theta$; since $[G,G]$ is central in $G$, $\Theta$ 
is a normal subgroup. Let $G'=G/\Theta$. Since $\Theta$ is discrete, the 
automorphisms of $G'$ are precisely factors on $G'$ of automorphisms $\tau$ of $G$
such that $\tau (\Theta) =\Theta$. We note that the shear automorphisms, 
namely those from $\Psi$  in the notation as in $\S\,$\ref{2step}, satisfy the 
condition. Therefore by Corollary~\ref{coraut} Aut$\,(G')$ is the 
semidirect product of $\Psi$ with the subgroup $\{\tau \in \tilde Z \tilde H \mid 
\tau (\Theta)=\Theta\}$. Now let $\tau =zh \in \tilde Z \tilde H $, where
$z \in \tilde Z $ and $h\in \tilde H $, be such that $\tau (\Theta)=\Theta$. 
Since $z$ acts by scalar multiplication on $[G,G]$, this implies that $\theta$ is
an eigenvector of $h$. Since $\theta \notin E$ it follows that $h$ is the identity 
element. Hence $\tau =z \in \tilde Z$, and since $\tau (\Theta)=\Theta$ it follows 
that the corresponding scalar transformation of $V=G/[G, G]$ is $\pm I$, where $I$ 
denotes the identity transformation. This implies 
the first assertion in the Corollary. The second one follows from this 
together with the fact that automorphisms from $\Psi$ factor to $G'$.  

\section{Anosov automorphisms}

As mentioned in the introduction study of the automorphism group can be applied
also to the general question of understanding the class of compact nilmanifolds
supporting Anosov automorphisms; see \cite{AuS}, \cite{D1}, \cite{De}. In this 
respect we note the following consequence of Theorem~\ref{thm}. 

\begin{corollary}\label{coran}
Let $G$ be the 2-step simply
connected nilpotent Lie group as in Theorem~\ref{thm} and $V=G/[G,G]$
(as before, realised as a vector space). If $\tau \in \hbox{Aut}\,(G)$ and
the factor $\overline \tau$ on $V$ has determinant 1 then $\overline \tau$ has 
a nonzero fixed point on $V$. Consequently, if  $\Gamma$ is a
lattice (a discrete cocompact subgroup) in $G$ and $\tau \in \hbox{Aut}\,(G)$ 
is such that $\tau (\Gamma)=\Gamma$ then the factor automorphism 
$\pi (\tau):G/\Gamma \to G/\Gamma$  is not an Anosov automorphism, and 
furthermore it is not ergodic. 
\end{corollary}

\noindent{\it Proof}:  We follow the notation as before. By Theorem~\ref{thm}
$\overline \tau \in ZH$, and furthermore if its determinant is 1 then 
$\overline \tau \in H$. The first assertion therefore follows from the fact
that for the $H$-action on $V$ every element of $h$ has a nonzero fixed point
in $V$; for hyperbolic elements this follows from consideration of weights,
for parabolic elements by unipotence, and for elliptic elements by odd-dimensionality 
of $V$ (analogous assertion holds for any irreducible representation of $H$ 
on an odd-dimensional vector space). 

If $\Gamma$ is a lattice in $G$ then 
$[G,G]\Gamma$ is closed and $[G,G]\Gamma/[G,G]$ is a lattice in $V=
G/[G,G]$; see \cite{R}. Therefore for any $\tau \in \hbox{Aut}\,(G)$ 
such that $\tau (\Gamma)=\Gamma$,  $\overline \tau$ has determinant~$\pm 1$.
Hence by the first part $\overline \tau^2$ has nonzero fixed points in $V$.
This implies that the factor of $\tau$ on $G/\Gamma$ is not an 
Anosov automorphism, and also that it is not ergodic. This proves the Corollary.

Compact nilmanifolds covered by the 3-step simply connected nilpotent
 Lie group 
$G$ as in $\S\,$\ref{nilauto} also do not support Anosov automorphisms; 
this is immediate from the fact that  there are no hyperbolic automorphisms, as 
$[G,[G,G]]$ is one-dimensional. 

\medskip
\noindent{\it Acknowledgement}: The author would like to thank Karel Dekimpe 
for useful comments on an earlier version of this manuscript.

\begin{flushleft}

\vskip2mm
Erwin Schr\"odinger Institute, Boltzmanngasse 9, A-1090 Vienna, Austria.

\vskip2mm 
{\it Permanent Address}: 
\vskip1pt 
School of Mathematics, Tata Institute of Fundamental Research, 
Homi Bhabha Road, Colaba, Mumbai 400 005, India\\ 
 
\vskip3pt 
{\tt E-mail: dani@math.tifr.res.in} 
\end{flushleft} 
\end{document}